\documentclass[11pt,a4paper]{article}

\usepackage{amsmath}
\usepackage{amsthm}
\usepackage{amssymb}
\usepackage{amscd}
\usepackage[all]{xy}

\title{Counterexamples to the Kawamata-Viehweg Vanishing \\
on Ruled Surfaces in Positive Characteristic
\footnote{This paper was partially supported by the National Natural Science 
Foundation of China (Grant No.\ 10901037) and Ph.D.\ Programs Foundation of 
Ministry of Education of China (Grant No.\ 20090071120004).}}
\author{Qihong Xie}
\date{}
\theoremstyle{plain}
\newtheorem{prop}{Proposition}[section]
\newtheorem{lem}[prop]{Lemma}
\newtheorem{thm}[prop]{Theorem}
\newtheorem{cor}[prop]{Corollary}
\newtheorem{conj}[prop]{Conjecture}
\newtheorem{prob}[prop]{Problem}
\newtheorem*{case}{Case}
\newtheorem*{subcase}{Subcase}

\theoremstyle{definition}
\newtheorem{defn}[prop]{Definition}
\newtheorem*{ack}{Acknowledgments}
\newtheorem*{conv}{Convention and Notation}

\theoremstyle{remark}
\newtheorem{rem}[prop]{Remark}
\newtheorem{ex}[prop]{Example}

\newcommand{\Q}{\mathbb Q}

\newcommand{\Z}{\mathbb Z}
\newcommand{\N}{\mathbb N}

\newcommand{\PP}{\mathbb P}
\newcommand{\OO}{\mathcal O}

\newcommand{\EE}{\mathcal E}
\newcommand{\NN}{\mathcal N}
\newcommand{\MM}{\mathcal M}
\newcommand{\LL}{\mathcal L}
\newcommand{\GG}{\mathcal G}

\newcommand{\BB}{\mathcal B}

\newcommand{\cS}{\mathcal S}
\newcommand{\cQ}{\mathcal Q}

\newcommand{\Supp}{\mathop{\rm Supp}\nolimits}
\newcommand{\Exc}{\mathop{\rm Exc}\nolimits}
\newcommand{\ch}{\mathop{\rm char}\nolimits}

\newcommand{\Spec}{\mathop{\rm Spec}\nolimits}
\newcommand{\Gal}{\mathop{\rm Gal}\nolimits}
\newcommand{\im}{\mathop{\rm im}\nolimits}
\newcommand{\ra}{\rightarrow}

\newcommand{\wt}{\widetilde}

\begin{document}
\maketitle

\begin{abstract}
We give counterexamples to the Kawamata-Viehweg vanishing theorem 
on ruled surfaces in positive characteristic, and 
prove that if there is a counterexample to the Kawamata-Viehweg 
vanishing theorem on a geometrically ruled surface $f:X\ra C$, then 
either $C$ is a Tango curve or all of sections of $f$ are ample.
\end{abstract}

\setcounter{section}{0}
\section{Introduction}\label{S1}

The study of pathology of algebraic geometry in positive 
characteristic is of certain interest, since pathology 
reveals some completely different geometric phenomena from 
those in complex geometry. 
For the celebrated Kodaira vanishing theorem, Raynaud \cite{ra} 
has given its counterexamples on quasi-elliptic surfaces and 
on surfaces of general type, which are smooth over an algebraically 
closed field $k$ of positive characteristic. Conversely, 
Tango \cite{ta72b} and Mukai \cite{mu79} have shown that if there is a 
counterexample to the Kodaira vanishing theorem on a smooth projective 
surface $X$, then $X$ must be a quasi-elliptic surface or a surface of 
general type. It turns out that if we restrict our attention to the 
category of smooth projective surfaces, then those well-known pathologies 
are very restrictive from the point of view of the classification theory 
of surfaces.

Mumford \cite{mu67} has shown that if the morphism 
$F^*:H^1(\OO_Y)\ra H^1(\OO_Y)$ is not injective, where $F^*$ is induced 
by the Frobenius map $F:Y\ra Y$ for a normal projective surface $Y$, 
then there is a normal projective surface $X$, a finite surjective 
separable morphism $\pi:X\ra Y$, and an ample line bundle $\LL$ on $X$ 
such that $H^1(X,\LL^{-1})\neq 0$. Therefore, in the category of normal 
projective surfaces, those pathologies exist widely.

For algebraic surfaces, it suffices to run the minimal model program in 
the category of smooth projective surfaces, while for higher-dimensional 
varieties, it is inevitable to consider either varieties with terminal 
singularities or log pairs with Kawamata log terminal singularities, 
where these singularities indeed arise naturally in the minimal model 
program. In the classification theory of three-dimensional varieties, 
we often take log pairs into account and use the adjunction theory 
to reduce the three-dimensional case to the surface case, thus it is 
necessary to investigate the birational structure of either singular 
surfaces or log surfaces instead of smooth surfaces.

As is well known, the Kawamata-Viehweg vanishing thereom generalizes 
the Kodaira vanishing theorem, and plays an essential role in the 
higher-dimensional minimal model program over the complex number field. 
Therefore, for the classification of three-dimensional varieties in 
positive characteristic, it is also necessary and interesting to find 
counterexamples to the Kawamata-Viehweg vanishing theorem on either 
singular surfaces or log surfaces, and give certain characterizations 
of these counterexamples.

A smooth projective curve $C$ is called a Tango curve if its Tango 
invariant $n(C)>0$, where the Tango invariant $n(C)$ was first introduced 
by Tango \cite{ta72a}. An important geometric property of a Tango curve 
$C$ is described as follows: there exists a rank two locally free sheaf 
$\EE$ on $C$ such that any quotient line bundle of $\EE$ has positive 
degree, however $\EE$ is not ample.

Note that the Kodaira vanishing theorem does hold on smooth ruled surfaces 
in positive characteristic (cf.\ \cite{ta72b}). However, it is shown that 
if $C$ is a Tango curve, then there exist counterexamples to the 
Kawamata-Viehweg vanishing theorem on a geometrically ruled surface $X$ 
over $C$ (cf.\ \cite[Example 3.7]{xie}).

The following is the main theorem in this paper, which is almost 
the converse of the above result (see Theorem \ref{2.2} for the 
Kawamata-Viehweg vanishing theorem).

\begin{thm}\label{1.1}
If there is a counterexample to the Kawamata-Viehweg vanishing theorem 
on a geometrically ruled surface $f:X\ra C$, then either $C$ is a 
Tango curve or all of sections of $f$ are ample.
\end{thm}

In fact, we may put forward the following conjecture, which would 
give a characterization of counterexamples to the Kawamata-Viehweg 
vanishing theorem on surfaces in positive characteristic.

\begin{conj}\label{1.2}
If there is a counterexample to the Kawamata-Viehweg vanishing theorem 
on a normal projective surface $X$, then there exists a dominant rational 
map $f:X\dashrightarrow C$ to a smooth projective curve $C$ such that 
$C$ is a Tango curve.
\end{conj}

The main idea of the proof of Theorem \ref{1.1} is to observe the 
behavior of cohomology classes of certain line bundles under the 
Frobenius map. It turns out that the construction of 
\cite[Example 3.7]{xie} (see also Theorem \ref{3.1}) and the proof 
of Theorem \ref{1.1} are inverse to each other to some extent. 
As a consequence of the main theorem, the Kawamata-Viehweg vanishing 
theorem holds on Hirzebruch surfaces. More generally, we can prove that 
the Kawamata-Viehweg vanishing theorem holds on rational surfaces 
\cite[Theorem 1.4]{xie07}, hence on log del Pezzo surfaces, 
which is an evidence for Conjecture \ref{1.2}.

In this paper, we also generalize slightly the construction of 
\cite[Example 3.7]{xie} to obtain counterexamples to the 
Kawamata-Viehweg vanishing theorem on certain general ruled surfaces 
and on certain normal projective surfaces of Picard number one with 
a nonrational singularity (see Corollaries \ref{3.3} and \ref{3.4}). 
All of results indicate that the pathologies of log surfaces do exist 
widely.

In \S \ref{S2}, we shall recall some definitions and elementary results. 
\S \ref{S3} may be regarded as a survey of counterexamples 
to the Kawamata-Viehweg vanishing theorem on ruled surfaces. 
\S \ref{S4} is devoted to the proof of the main theorem. 

\begin{conv}
Throughout this paper, 
{\it we always work over an algebraically closed field $k$ 
of characteristic $p>0$} unless otherwise stated. 
A surface $X$ is called a ruled surface over a smooth curve $C$ if 
there is a morphism $f:X\ra C$ such that $f_*\OO_X=\OO_C$ and the 
general fibre of $f$ is $\PP^1$, and it is called a geometrically 
ruled surface if $X$ is a $\PP^1$-bundle over $C$. 
We use $\equiv$ to denote numerical equivalence, 
$[B]=\sum [b_i] B_i$ to denote the round-down of 
a $\Q$-divisor $B=\sum b_iB_i$, and $K(C)$ to 
denote the rational function field of a curve $C$.
\end{conv}

\begin{ack}
I would like to express my gratitude to Professors Yujiro Kawamata 
and Takao Fujita for valuable advice and warm encouragement. 
I would also like to thank Professors Stefan Schr\"oer and Shigeru 
Mukai for useful suggestions. I am very grateful to the referee for 
carefully reading the manuscript and giving many useful comments, 
which make this paper more readable.
\end{ack}

\section{Preliminaries}\label{S2}

\begin{defn}\label{2.1}
Let $X$ be a normal proper algebraic variety over an algebraically 
closed field $k$, and $B=\sum b_iB_i$ an effective $\Q$-divisor on $X$. 
The pair $(X,B)$ is said to be Kawamata log terminal (KLT, for short), 
or to have Kawamata log terminal singularities, 
if the following conditions hold:

(i) $K_X+B$ is $\Q$-Cartier, i.e.\ $r(K_X+B)$ is Cartier for some $r\in\N$;

(ii) For any birational morphism $f:Y\rightarrow X$, if we write 
$K_Y+B_Y\equiv f^*(K_X+B)$, where $B_Y=\sum a_iE_i$ is a $\Q$-divisor on $Y$, 
then $a_i<1$ hold for all $i$.
\end{defn}

First, it follows from (ii) that $[B]=0$, i.e.\ $b_i<1$ for all $i$. 
Secondly, this definition is characteristic free. 
Thirdly, when $\ch(k)=0$ or $\dim X\leq 2$, the pair $(X,B)$ admits a log 
resolution, i.e.\ there exists a desingularization $f:Y\rightarrow X$ 
from a nonsingular variety $Y$, such that the union of the strict 
transform $f_*^{-1}B$ of $B$ and the exceptional locus $\Exc(f)$ of $f$ 
has simple normal crossing support. We can show that when $\ch(k)=0$ or 
$\dim X\leq 2$, (ii) holds for all birational morphisms is equivalent to 
that (ii) holds for a log resolution of $(X,B)$.

Let us mention a simple example of KLT pair. Let $X$ be a nonsingular 
variety and $B$ an effective $\Q$-divisor on $X$ such that 
$[B]=0$ and $\Supp(B)$ is simple normal crossing. 
Then $(X,B)$ is KLT (see \cite[Corollary 2.31]{km} for the proof).

The Kawamata-Viehweg vanishing theorem is one of the most 
important generalizations of the Kodaira vanishing theorem 
(cf.\ \cite[Theorem 1-2-5]{kmm}).

\begin{thm}[Kawamata-Viehweg vanishing]\label{2.2}
Let $X$ be a normal projective algebraic variety over an algebraically 
closed field $k$ with $\ch(k)=0$, $B=\sum b_iB_i$ an effective 
$\Q$-divisor on $X$, and $D$ a $\Q$-Cartier Weil divisor on $X$. 
Assume that $(X,B)$ is KLT and that $H=D-(K_X+B)$ is ample. 
Then $H^i(X,D)=0$ holds for any $i>0$.
\end{thm}

From now on, we always assume that $C$ is a smooth projective 
curve over an algebraically closed field $k$ of characteristic $p>0$. 
First, we recall the following definition from 
\cite[Definitions 9 and 11]{ta72a}.

\begin{defn}\label{2.3}
Let $f\in K(C)$ be a rational function on $C$. 
\[ n(f):=\deg\biggl{[}\frac{(df)}{p}\bigg{]}, \]
where $(df)=\sum_{x\in C}v_x(df)x$ is the divisor associated to 
the rational differential 1-form $df$. The Tango invariant 
$n(C)$ of the curve $C$ is defined by
\[ n(C):=\max\{n(f) \, | \, f\in K(C),f\not\in K^p(C) \}. \]
\end{defn}

If $f\not\in K^p(C)$, then $(df)$ is a canonical divisor on $C$ 
with degree $2(g(C)-1)$. It is easy to see that $n(C)\leq [2(g(C)-1)/p]$. 
Therefore, if $n(C)>0$ then $g(C)\geq 2$ holds. On the other hand, 
it follows from \cite[Proposition 14]{ta72a} that $n(C)\geq 0$ holds if 
$g(C)\geq 1$. By an elementary calculation, we can show that $n(\PP^1)=-1$ 
and $n(C)=0$ holds for any elliptic curve $C$.

\begin{ex}\label{2.4}
There do exist smooth projective curves $C$ such that $n(C)>0$ 
for each characteristic $p>0$.

(i) Assume $p\geq 3$. Let $h\geq 3$ be an odd integer, and let $C$ be 
the projective completion at infinity of the affine curve defined by 
$y^2=x^{ph}+x^{p+1}+1$. It is easy to verify that $C$ is a smooth 
hyperelliptic curve and that $(d(y/x^p))=(ph-3)z_\infty$, 
where $z_\infty$ is the infinity point of $C$ 
(cf.\ \cite[Ch.\ III, \S 6.5]{sha}). Hence $n(C)=n(y/x^p)=h-1>0$.

(ii) Let $h>2$ be an integer, and let $C$ be the projective completion 
at infinity of the Artin-Schreier cover of the affine line 
defined by $y^{hp-1}=x^p-x$. It is easy to verify that $C$ 
is a smooth curve and that $(dy)=p(h(p-1)-2)z_\infty$, 
where $z_\infty$ is the infinity point of $C$. 
Hence $n(C)=n(y)=h(p-1)-2>0$ (cf.\ \cite{ra}).

(iii) Assume $p\geq 3$. Let $C\subset\PP^2$ be the projective curve 
defined by $x_0^{p+1}=x_1x_2(x_0^{p-1}+x_1^{p-1}-x_2^{p-1})$. 
Then we can show that $C$ is a smooth curve with $n(C)=n(x_0/x_1)=p-2>0$ 
(cf.\ \cite[Proposition 28]{ta72a}).
\end{ex}

Let $F:C\rightarrow C$ be the Frobenius map. We have the following 
exact sequences of locally free sheaves on $C$ 
(cf.\ \cite[Proposition 3]{ta72a}):
\begin{eqnarray}
0\rightarrow \OO_C\rightarrow F_*\OO_C\rightarrow \BB^1\rightarrow 0, 
\label{es:1} \\
0\rightarrow \BB^1\rightarrow F_*\Omega^1_C\stackrel{c}{\rightarrow} 
\Omega^1_C\rightarrow 0, \label{es:2} 
\end{eqnarray}
where $\BB^1$ is the image of the map $F_*(d):F_*\OO_C\rightarrow 
F_*\Omega^1_C$, and $c$ is the Cartier operator.

\begin{lem}\label{2.5}
Notation as above, let $L$ be a divisor on $C$. 
Then $H^0(C,\BB^1(-L))=\{ df \,|\, f\in K(C),\,(df)\geq pL \}$. 
In particular, $n(C)>0$ if and only if there is an ample divisor 
$L$ on $C$ such that $H^0(C,\BB^1(-L))\neq 0$.
\end{lem}

\begin{proof}
Tensoring (\ref{es:2}) by $\OO_C(-L)$, we have:
\[ 0\rightarrow \BB^1(-L)\rightarrow F_*(\Omega^1_C(-pL))
\stackrel{c(-L)}{\rightarrow} \Omega^1_C(-L)\rightarrow 0. \]
By the definition of $\BB^1$ and $H^0(C,\Omega^1_C(-pL))=
\{ \omega\in\Omega^1_C \,|\, (\omega)\geq pL \}$, 
we have $H^0(C,\BB^1(-L))=\{ df \,|\, f\in K(C),\,(df)\geq pL \}$.

If $n(C)>0$, then there exists an $f_0\in K(C)$ such that 
$n(f_0)=\deg [(df_0)/p]=n(C)>0$. Let $L=[(df_0)/p]$. Then we have 
$\deg L=n(C)>0$ and $(df_0)\geq pL$, hence $0\neq df_0\in H^0(C,\BB^1(-L))$.

If there is an ample divisor $L$ on $C$ such that $H^0(C,\BB^1(-L))\neq 0$, 
then we can take an $f_0\in K(C)$ with $(df_0)\geq pL$, thus we have 
$n(f_0)\geq\deg L>0$, hence $n(C)>0$.
\end{proof}

\begin{defn}\label{2.6}
A smooth projective curve $C$ is called a Tango curve if $n(C)>0$ holds.
Let $C$ be a Tango curve and $f_0\in K(C)$ such that $n(f_0)=n(C)$. 
Then the divisor $L=[(df_0)/p]$ is called a base divisor of $C$. 
A base divisor $L$ is said to be divisible if $\frac{1}{2}L$ is integral 
when $p\geq 3$, or $\frac{1}{3}L$ is integral when $p=2$.
A Tango curve $C$ is said to be of integral type if there is a 
divisible base divisor $L$ of $C$.

A smooth projective curve $C$ is called a Raynaud-Tango curve, 
if there is an $f_0\in K(C)$ and an integral divisor $L$ on $C$ 
such that $(df_0)=pL$ and $\deg L>0$. If this is the case, then 
$K_C=pL$, $n(C)=n(f_0)=\deg L>0$ and $L$ is a base divisor of $C$.

For instance, the curve $C$ given in Example \ref{2.4} (i) is a 
Tango curve of integral type, and the curve $C$ given in Example 
\ref{2.4} (ii) is a Raynaud-Tango curve of integral type for some 
suitable $h$ when $p=2$, while the curve $C$ given in Example 
\ref{2.4} (iii) is a Tango curve, which is not of integral type.
\end{defn}

\section{Counterexamples to the Kawamata-Viehweg vanishing}\label{S3}

In \cite[Example 3.7]{xie}, some counterexamples to the Kawamata-Viehweg 
vanishing theorem on certain geometrically ruled surfaces have been 
constructed. We recall and simplify the explicit construction for 
the convenience of the reader.

\begin{thm}\label{3.1}
If $C$ is a Tango curve, then there is a $\PP^1$-bundle 
$f:X\rightarrow C$, an effective $\Q$-divisor $B$ and an 
integral divisor $D$ on $X$ such that $(X,B)$ is KLT and 
$H=D-(K_X+B)$ is ample. However, we have $H^1(X,D)\neq 0$.
\end{thm}

\begin{proof}
By Lemma \ref{2.5}, we can take a base divisor $L$ of $C$ 
such that $\deg L=n(C)>0$ and $H^0(C,\BB^1(-L))\neq 0$. 
Let $\LL=\OO_C(L)$. Tensoring the exact sequence (\ref{es:1}) 
by $\LL^{-1}$ and taking cohomology groups, 
we have the following exact sequence:
\[ 0\rightarrow H^0(C,\BB^1(-L))\stackrel{\eta}{\rightarrow} 
H^1(C,\LL^{-1})\stackrel{F^*}{\rightarrow} H^1(C,\LL^{-p}). \]
Since $\eta$ is injective, we can take an element $0\neq \alpha\in 
H^0(C,\BB^1(-L))$ such that $0\neq\eta(\alpha)\in H^1(C,\LL^{-1})$, 
which determines the following extension of $\LL$ by $\OO_C$: 
\begin{eqnarray}
0\rightarrow \OO_C\rightarrow \EE\rightarrow \LL\rightarrow 0. \label{es:3}
\end{eqnarray}
Pulling back the exact sequence (\ref{es:3}) by the Frobenius map $F$, 
we have the following split exact sequence:
\begin{eqnarray}
0\rightarrow \OO_C\rightarrow F^*\EE\rightarrow \LL^p\rightarrow 0, 
\label{es:4}
\end{eqnarray}
since the extension class of (\ref{es:4}) is just $F^*\eta(\alpha)=0$.

Let $X=\PP(\EE)$ be the $\PP^1$-bundle over $C$, $f:X\rightarrow C$ 
the projection, $\OO_X(1)$ the tautological line bundle and $F$ the 
fibre of $f$. The sequence (\ref{es:3}) determines a section $E$ of $f$ 
through the map $H^0(C,\OO_C)\hookrightarrow H^0(C,\EE)=H^0(X,\OO_X(1))$, 
hence $\OO_X(E)\cong \OO_X(1)$. 
The sequence (\ref{es:4}) induces an exact sequence:
\[ 0\rightarrow \OO_C\rightarrow F^*\EE\otimes\LL^{-p}\rightarrow 
\LL^{-p}\rightarrow 0, \]
which determines a section $t\in H^0(X,\OO_X(p)\otimes f^*\LL^{-p})$ 
through the maps $H^0(C,\OO_C)\hookrightarrow H^0(C,F^*\EE\otimes\LL^{-p})
\hookrightarrow H^0(C,S^p(\EE)\otimes\LL^{-p})=
H^0(X,\OO_X(p)\otimes f^*\LL^{-p})$. The section $t$ determines an 
irreducible curve $C'$ on $X$ such that $\OO_X(C')\cong 
\OO_X(p)\otimes f^*\LL^{-p}$. It is easy to verify that 
both $E$ and $C'$ are smooth over $k$, and $E\cap C'=\emptyset$.

Let $c$ be a rational number such that $1/p<c<1$ and $cp\not\in\Z$. Then 
$q=[cp]-1$ is a nonnegative integer. Let $B=cC'$ and $D=qE+f^*(K_C-qL)$. 
Then we have $H=D-(K_X+B)\equiv ([cp]+1-cp)E+(cp-[cp])f^*L$. 
It is easy to see that $(X,B)$ is KLT. Since $E^2=\deg\EE=\deg\LL>0$, 
$E$ is a nef and big divisor on $X$. 
Furthermore, since $E$ is $f$-ample, $H$ is an ample $\Q$-divisor on $X$.

Next we show $H^1(X,D)\neq 0$. Since $D|_F$ is nef on $F\cong\PP^1$, 
we have $R^1f_*\OO_X(D)=0$ and $f_*\OO_X(D)$ is a locally free sheaf 
on $C$. By the Leray spectral sequence, we have 
\begin{eqnarray}
H^1(X,D) & = & H^1(C,f_*\OO_X(D)) \nonumber \\
& = & H^0(C,(f_*\OO_X(D))^\vee\otimes\omega_C)^\vee
=H^0(C,(f_*\OO_X(D-f^*K_C))^\vee)^\vee \nonumber \\
& = & H^0(C,S^q(\EE)^\vee\otimes\LL^q)^\vee. \nonumber
\end{eqnarray}
Since $\LL^q$ is a quotient sheaf of $S^q(\EE)$, 
$\LL^{-q}$ is a subsheaf of $S^q(\EE)^\vee$. 
Thus we have
\[ H^1(X,D)^\vee=H^0(C,S^q(\EE)^\vee\otimes\LL^q)
\supseteq H^0(C,\LL^{-q}\otimes\LL^q)=H^0(C,\OO_C)=k, \]
which is desired.

Finally, we give some common choices of $c,q,B$ and $D$ for later use:
\begin{enumerate}
\item[(\dag 1)] $p\geq 3$: $c=1/2$, $q=\frac{p-3}{2}$, 
$B=\frac{1}{2}C'$, $D=\frac{p-3}{2}E+f^*(K_C-\frac{p-3}{2}L)$.

\item[(\dag 2)] $p=2$: $c=2/3$, $q=0$, $B=\frac{2}{3}C'$, $D=f^*K_C$.

\item[(\dag 3)] $p=3$: $c=5/6$, $q=1$, $B=\frac{5}{6}C'$, $D=E+f^*(K_C-L)$.
\end{enumerate}
\end{proof}

Note that in Theorem \ref{3.1}, $X$ is a geometrically ruled surface 
and $\Supp(B)$ is a smooth curve. We can generalize this counterexample 
slightly to a new one where $X$ is a general ruled surface and $\Supp(B)$ 
is simple normal crossing. First, we need the following vanishing result 
\cite[Corollary 2.2.5]{kk}, which holds in arbitrary characteristic.

\begin{lem}\label{3.2}
Let $h: \wt{X}\ra X$ be a proper birational morphism between normal surfaces 
with $\wt{X}$ smooth and with exceptional locus $E=\cup_{i=1}^s E_i$. 
Let $L$ be an integral divisor on $\wt{X}$, $0\leq b_1,\cdots,b_s<1$ 
rational numbers, and $N$ an $h$-nef $\Q$-divisor on $\wt{X}$. Assume 
$ L\equiv K_{\wt{X}}+\sum_{i=1}^s b_iE_i+N$. Then we have 
$R^1h_*\OO_{\wt{X}}(L)=0$.
\end{lem}

\begin{cor}\label{3.3}
With the same notation and construction as in Theorem \ref{3.1}, 
let $h:\wt{X}\ra X$ be a proper birational morphism such that $\wt{X}$ 
is smooth and $h_*^{-1}B\cup\Exc(h)$ has simple normal crossing support, 
then there exists a counterexample to the Kawamata-Viehweg vanishing 
theorem on $\wt{X}$.
\end{cor}

\begin{proof}
Let $\wt{f}=f\circ h: \wt{X}\rightarrow C$ be the induced 
morphism and write 
\[ K_{\wt{X}}+h_*^{-1}B \equiv h^*(K_X+B)+\sum_{i=1}^s a_iE_i, \]
where $E_i$ are exceptional curves of $h$ and $a_i>-1$ for all 
$1\leq i\leq s$.

Denote $\ulcorner a\urcorner$ (resp.\ $\{ a\}$) the round up (resp.\ 
fractional part) of a real number $a$. 
Let $\wt{D}=h^*D+\sum_{i=1}^s \ulcorner a_i \urcorner E_i$. Then we have 
\[ \wt{D} \equiv K_{\wt{X}}+\sum_{i=1}^s\{ -a_i\}E_i+h^*(D-(K_X+B))
+h_*^{-1}B. \]
Since $h^*(D-(K_X+B))+h_*^{-1}B$ is $h$-nef, 
we have $R^1h_*\OO_{\wt{X}}(\wt{D})=0$ by Lemma \ref{3.2}.

We can take $0<\delta_1,\cdots,\delta_s\ll 1$ such that 
\begin{enumerate}
\item[(i)] $\{ -a_i\}+\delta_i<1$ holds for any $1\leq i\leq s$, and 

\item[(ii)] $h^*(D-(K_X+B))-\sum_{i=1}^s\delta_iE_i$ is ample.
\end{enumerate}

Let $\wt{B}=h_*^{-1}B+\sum_{i=1}^s(\{ -a_i\}+\delta_i)E_i$. 
Then $(\wt{X},\wt{B})$ is KLT and $\wt{D}-(K_{\wt{X}}+\wt{B})\equiv 
h^*(D-(K_X+B))-\sum_{i=1}^s\delta_iE_i$ is ample. 
By the Leray spectral sequence and the projection formula, 
we have $H^1(\wt{X},\wt{D})\cong H^1(X,h_*\OO_{\wt{X}}(\wt{D}))
=H^1(X,D)\neq 0$. Thus we obtain a counterexample to the 
Kawamata-Viehweg vanishing theorem on $\wt{X}$.
\end{proof}

\begin{rem}\label{3.6}
In a similar manner, we can generalize \cite[Examples 3.9 and 3.10]{xie} 
slightly to new ones which are counterexamples to the logarithmic Koll\'{a}r 
vanishing and the logarithmic semipositivity theorems on certain general 
ruled surfaces.
\end{rem}

From the construction in Theorem \ref{3.1}, we can obtain a normal 
projective surface of Picard number one with a nonrational singularity, 
on which the Kawamata-Viehweg vanishing theorem fails.

\begin{cor}\label{3.4}
Let $C$ be a Raynaud-Tango curve with $\ch(k)=p\geq 3$. 
Then there is a normal projective surface $Y$ of Picard number one 
with a nonrational singularity, and a $\Q$-Cartier Weil divisor $D'$ 
on $Y$ such that $D'-K_Y$ is ample. However, we have $H^1(Y,D')\neq 0$.
\end{cor}

\begin{proof}
We follow the same notation and construction as in Theorem \ref{3.1}. 
Since $E$ is irreducible and $E^2>0$, by \cite[Theorem 1.10]{fu}, 
$E$ is semiample, i.e.\ a suitable multiple of $E$ is 
base-point-free and defines a birational morphism $\pi: X\ra Y$ 
such that $Y$ is a normal projective surface. As a consequence, 
$\pi_*E$ is an ample $\Q$-Cartier Weil divisor on $Y$. 
Since $\rho(X)=2$, the exceptional locus of $\pi$ consists of only 
one irreducible curve, which is just $C'$ since $E\cap C'=\emptyset$. 
Hence the Picard number $\rho(Y)=1$. Let $y=\pi_*(C')\in Y$. 
Then $y$ is the only one singular point of $Y$, which is nonrational 
since $g(C')=g(C)\geq 2$.

Let $n=n(C)=\deg L$ and $g=g(C)$. Since $C$ is a Raynaud-Tango curve, 
we have $K_C=pL$ and $2g-2=pn$. 
The following constants $a,b$ are introduced to simplify the argument.
\begin{itemize}
\item When $p\geq 5$, we define $a=(p+1)/p$ and $b=(p+3)/2p$ and we use 
the choice (\dag 1) in Theorem \ref{3.1}. In particular, we have 
$D=\frac{p-3}{2}E+f^*(K_C-\frac{p-3}{2}L)$.

\item When $p=3$, we define $a=4/3$ and $b=2/3$ and we use the choice 
(\dag 3) in Theorem \ref{3.1}. In particular, we have $D=E+f^*(K_C-L)$.
\end{itemize}

By using the relations $K_X=-2E+f^*(K_C+L)$ and $C'=pE-pf^*L$ 
and through a direct calculation in the corresponding choices, 
we can easily show that:
\begin{enumerate}
\item[(i)] $p(K_X+aC')$ is linearly equivalent to $p(p-1)E$, hence 
$K_Y=\pi_*(K_X+aC')$ is ample $\Q$-Cartier and we have 
$K_X\equiv \pi^*K_Y-aC'$;

\item[(ii)] $p(D+bC')$ is linearly equivalent to $p^2E$, hence 
$\pi_*D=\pi_*(D+bC')$ is $\Q$-Cartier and we have 
$D\equiv \pi^*\pi_*D-bC'$.
\end{enumerate}

It follows easily from (i) that $Y$ is not KLT. Define $D'=\pi_*D$. 
Then $D'$ is a $\Q$-Cartier Weil divisor on $Y$ by (ii). 
Note that $D'-K_Y=\pi_*(D-K_X)=\pi_*(H+B)=\pi_*H$ is ample on $Y$.

Since $D\equiv \pi^*\pi_*D-bC'$ with $0<b<1$, we have $[\pi^*\pi_*D]=D$. 
By the projection formula, we have $\pi_*\OO_X(D)=\OO_Y(\pi_*D)=\OO_Y(D')$. 
Since $D=K_X+cC'+H$, where $0<c<1$ and $H$ is ample, 
we have $R^1\pi_*\OO_X(D)=0$ by Lemma \ref{3.2}. 
By the Leray spectral sequence, 
$H^1(Y,D')=H^1(Y,\pi_*\OO_X(D))=H^1(X,D)\neq 0$.
\end{proof}

From the construction in Theorem \ref{3.1}, we can give a more general 
statement which generalizes and clarifies, to some extent, the construction 
of counterexamples to the Kodaira vanishing theorem given in \cite{ra}. 
We sketch out the proof for the convenience of the reader.

\begin{thm}\label{3.5}
If $C$ is a Tango curve of integral type, then there is a fibration 
$g:Y\ra C$ from a smooth projective surface $Y$, such that 
there are counterexamples to the Kodaira vanishing, 
the Koll\'{a}r vanishing and the semipositivity theorems on $Y$.
\end{thm}

\begin{proof}
With the same notation and construction as in Theorem \ref{3.1}, 
we can obtain a $\PP^1$-bundle $f:X\ra C$ over the Tango curve $C$ 
of integral type. Let $L$ be a divisible base divisor of $C$, and let 
\[ N=
\begin{cases}
\frac{1}{2}L & \text{if $p\geq 3$}, \\
\frac{1}{3}L & \text{if $p=2$}.
\end{cases}
\]
Then by definition, $N$ is an integral divisor on $C$. Let 
\[
M=
\begin{cases}
\frac{1}{2}(p+1)E-pf^*N & \text{if $p\geq 3$}, \\
E-2f^*N & \text{if $p=2$}.
\end{cases}
\]
Then $M$ is an integral divisor on $X$ such that 
$\OO_X(M)^2\cong\OO_X(E+C')$ if $p\geq 3$ or 
$\OO_X(M)^3\cong\OO_X(E+C')$ if $p=2$. Thus we can take a cyclic cover 
$\varphi:Y\ra X$ of degree 2 or 3 associated to the line bundle $\OO_X(M)$ 
with branch locus $E\cup C'$, which is the same as in \cite{ra}. 
Let $g=f\circ\varphi: Y\ra C$ be the induced morphism. 
Since both $E$ and $C'$ are smooth over $k$ and $E\cap C'=\emptyset$, 
$Y$ is a smooth projective surface over $k$ by \cite[Lemma 3.15 (b)]{ev}.

Let $\wt{E}=(\varphi^*E)_{\mathrm{red}}$ and $\wt{D}=K_Y+\wt{E}+g^*N$. 
Then through a direct calculation, we can verify that $\wt{D}-K_Y$ 
is ample and $H^1(Y,\wt{D})\neq 0$, which is a counterexample to 
the Kodaira vanishing theorem on $Y$.

Let $\NN=\OO_C(N)$. Then through a direct calculation, we can prove that
\[
\NN\otimes g_*\omega_Y=
\begin{cases}
S^{(p-3)/2}(\EE)\otimes \NN^{3-p}\otimes\omega_C & \text{if $p\geq 3$}, \\
\omega_C & \text{if $p=2$}, 
\end{cases}
\]
hence $H^1(C,\NN\otimes g_*\omega_Y)\neq 0$ holds, which is a 
counterexample to the Koll\'{a}r vanishing theorem for $g:Y\ra C$ 
(cf.\ \cite[Theorem 2.1]{ko}).

Finally, through a direct calculation, we can show that $g_*\omega_{Y/C}$ 
has a quotient line bundle $\NN^{-1}$ of negative degree, which 
is a counterexample to the semipositivity theorem for $g:Y\ra C$ 
(cf.\ \cite[Theorem 2.7]{fu78} and the general statement 
\cite[Theorem 5]{ka81}).
\end{proof}

We may regard \cite[Examples 3.7, 3.9 and 3.10]{xie} as logarithmic 
counterparts of Theorem \ref{3.5} in some sense. However, in general, 
it seems impossible to produce further counterexamples in the 
absolute case by means of blowing-ups and divisorial contractions 
from the constructions in Theorem \ref{3.5}.

\section{Proof of the main theorem}\label{S4}

In this section, we shall prove the main theorem in the following 
explicit form.

\begin{thm}\label{4.1}
Let $f:X\ra C$ be a $\PP^1$-bundle over a smooth 
projective curve $C$. Let $D$ be an integral divisor on $X$, 
$B$ an effective $\Q$-divisor on $X$ such that $(X,B)$ is KLT and 
that $H=D-(K_X+B)$ is ample. Assume $H^1(X,D)\neq 0$. Then either 

(i) $C$ is a Tango curve, or 

(ii) all of sections of $f$ are ample.
\end{thm}

In fact, the proof of Theorem \ref{4.1} actually gives more precise 
information on when the Kawamata-Viehweg vanishing theorem holds on 
a $\PP^1$-bundle over a curve, which is summarized in Corollary \ref{4.6}. 
Before proving Theorem \ref{4.1}, we need some auxiliary lemmas. 
The first one is a special case of a result in \cite{ta72b}, while 
we give a complete proof here for the convenience of the reader.

\begin{lem}\label{4.2}
Let $f:X\ra C$ be a $\PP^1$-bundle over a smooth projective curve $C$, 
and $\LL$ an ample line bundle on $X$. Then $H^1(X,\LL^{-1})=0$.
\end{lem}

\begin{proof}
Let $F:X\ra X$ be the Frobenius map. Then we have the following exact 
sequences of $\OO_X$-modules:
\[
\xymatrix{
0 \ar[r] & \OO_X \ar[r] & F_*\OO_X \ar[r] \ar[dr]_{F_*(d)} 
& \BB_X^1 \ar[r] \ar@{^{(}->}[d] & 0, \\
         &              &                                  
& F_*\Omega_X^1                  &      
}
\]
where $\BB_X^1=\im(F_*(d):F_*\OO_X\ra F_*\Omega_X^1)$.

By means of the above sequences, we have the following implications 
for all ample line bundles $\LL$:
$H^0(X,\Omega_X^1\otimes\LL^{-1})=0$
$\Longrightarrow$
$H^0(X,F_*\Omega_X^1\otimes\LL^{-1})=0$
$\Longrightarrow$
$H^0(X,\BB_X^1\otimes\LL^{-1})=0$
$\Longrightarrow$
$H^1(X,\LL^{-1})\ra H^1(X,\LL^{-p})$ is injective
$\Longrightarrow$
$H^1(X,\LL^{-1})=0$ by Serre vanishing. 
Next we show that $H^0(X,\Omega_X^1\otimes\LL^{-1})=0$ holds for any 
ample line bundle $\LL$ on $X$, which will complete the proof of Lemma 
\ref{4.2}.

Since $f:X\ra C$ is smooth, we have the following exact sequence:
\[
0\ra f^*\Omega^1_C\ra \Omega^1_X\ra \Omega^1_{X/C}\ra 0.
\]
Tensoring the above sequence by $\LL^{-1}$ and taking cohomology groups, 
it is easy to see that to prove $H^0(X,\Omega_X^1\otimes\LL^{-1})=0$, it 
suffices to prove $H^0(X,f^*\Omega^1_C\otimes\LL^{-1})=H^0(X,\Omega^1_{X/C}
\otimes\LL^{-1})=0$. Let $F$ be the fiber of $f$. Note that $f^*\Omega^1_C
\otimes\LL^{-1}$ and $\Omega^1_{X/C}\otimes\LL^{-1}$ are line bundles on $X$. 
If there is an effective divisor $D$ on $X$ such that $\OO_X(D)\cong 
f^*\Omega^1_C\otimes\LL^{-1}$, then we have $\deg_F f^*\Omega^1_C\otimes
\LL^{-1}=D.F\geq 0$, but $\deg_F f^*\Omega^1_C\otimes\LL^{-1}=\deg_F\LL^{-1}
<0$ since $\LL$ is ample, which is a contradiction. Therefore we have 
$H^0(X,f^*\Omega^1_C\otimes\LL^{-1})=0$. The verification of $H^0(X,
\Omega^1_{X/C}\otimes\LL^{-1})=0$ is similar.
\end{proof}

\begin{lem}\label{4.3}
Let $C$ be a smooth projective curve, and $\LL$ a line bundle on $C$ 
with $\deg\LL<0$. Let $\EE$ be a locally free sheaf on $C$ fitting into 
the following exact sequence:
\begin{eqnarray}
0\ra \OO_C\ra \EE\ra \LL\ra 0. \label{es:5}
\end{eqnarray}
Then for any $n\geq 0$ and any quotient line bundle $\GG$ of $S^n(\EE)$, 
we have $\deg\GG\geq n\deg\LL$.
\end{lem}

\begin{proof}
Consider the exact sequence $S^n(\EE)\ra \GG\ra 0$. If $n=0$, 
then we have $\GG\cong\OO_C$ and the conclusion is trivial. 
Assume $n\geq 1$. Let $\alpha\in H^1(C,\LL^{-1})$ be the extension 
class of the exact sequence (\ref{es:5}). After pulling back (\ref{es:5}) 
by the iterated Frobenius map $F^m$, we have the following exact sequence:
\begin{eqnarray}
0\ra \OO_C\ra F^{m*}\EE\ra \LL^{p^m}\ra 0. \label{es:6}
\end{eqnarray}
We claim that there exists an $m_0\in\N$ such that the exact sequence 
(\ref{es:6}) is split for any $m\geq m_0$. Indeed, the extension class 
of (\ref{es:6}) is just $F^{m*}(\alpha)\in H^1(C,\LL^{-p^m})$. Since 
$\LL^{-1}$ is ample, there exists an $m_0\in\N$ such that the cohomology 
group $H^1(C,\LL^{-p^m})=0$ for any $m\geq m_0$. 
Therefore, $F^{m*}(\EE\otimes\LL^{-1})\cong \LL^{-p^m}\oplus\OO_C$ 
for any $m\geq m_0$.

Let $\MM$ be an ample line bundle on $C$. Tensoring the exact sequence 
$S^n(\EE)\ra \GG\ra 0$ by $\LL^{-n}$, pulling back by the iterated 
Frobenius map $F^m$ and tensoring by $\MM^n$ successively, we have 
the following exact sequence:
\[ S^n(F^{m*}(\EE\otimes\LL^{-1})\otimes\MM)\ra 
(\GG\otimes\LL^{-n})^{p^m}\otimes\MM^n\ra 0. \]
Since $F^{m*}(\EE\otimes\LL^{-1})\otimes\MM\cong (\LL^{-p^m}\otimes\MM)
\oplus\MM$ is a direct sum of ample line bundles on $C$, 
$S^n(F^{m*}(\EE\otimes\LL^{-1})\otimes\MM)$ is also a direct sum of 
ample line bundles by \cite[Lemma 6.5]{ha66}. 
Hence $S^n(F^{m*}(\EE\otimes\LL^{-1})\otimes\MM)$ is ample 
by \cite[Proposition 2.2]{ha66}, and 
its quotient $(\GG\otimes\LL^{-n})^{p^m}\otimes\MM^n$ is also ample
by \cite[Proposition 2.2]{ha66} again. 
Thus $p^m(\deg\GG-n\deg\LL)+n\deg\MM>0$ holds for any $m\geq m_0$. 
As a consequence, we have $\deg\GG\geq n\deg\LL$.
\end{proof}

\begin{lem}\label{4.4}
Let $C$ be a smooth projective curve, and $f:X=\PP(\EE)\ra C$ a 
$\PP^1$-bundle associated to a rank two locally free sheaf $\EE$ on $C$. 
Let $\varphi:B\ra C$ be a surjective morphism from a smooth projective 
curve $B$, and $X_B=X\times_C B$ the fibre product. Let $R$ be an irreducible 
curve on $X$, and $R_B=R\times_C B$ the corresponding curve on $X_B$. 
Then the following assertions hold:

(i) $g:X_B\ra B$ is the $\PP^1$-bundle associated to the locally free sheaf 
$\varphi^*\EE$ on $B$.

(ii) If $R$ is a section of $f:X\ra C$, then $R_B$ is a section of 
$g:X_B\ra B$.

(iii) Assume that both $\varphi:B\ra C$ and $f|_R:R\ra C$ are purely 
inseparable of degree $p^n$. Then $R_B$ is a section of $g:X_B\ra B$. 
In particular, if $f|_R:R\ra C$ is of degree 1 then $R$ is a section of $f$.
\[
\xymatrix{
R_B\subset X_B \ar[d]_g \ar[r]^\psi & R\subset X \ar[d]^f \\
B \ar[r]^\varphi & C 
}
\]
\end{lem}

\begin{proof}
(i) It is obvious by the functoriality of projective space bundles.

(ii) If $R$ is a section of $f$, then $R_B=R\times_C B\cong B$, 
i.e.\ $R_B$ is a section of $g$.

(iii) $\varphi:B\ra C$ is purely inseparable of degree $p^n$, i.e.\ 
$K(B)/K(C)$ is a purely inseparable field extension of degree $p^n$. 
By \cite[Proposition IV.2.5]{ha}, we can assume that 
$\varphi:B=C_{p^n}\ra C$ is the iterated $k$-linear Frobenius map. 
The conclusion of (iii) can be verified locally, so we may assume that 
$C=\Spec A$, $\varphi$ is induced by the iterated $k$-linear Frobenius 
map $\varphi^*:A\ra A$, and $X=\Spec A[x]$, where $x$ is an indeterminant. 
Hence $\psi:X_B\ra X$ is induced by the ring homomorphism 
$\psi^*:A[x]\ra A[x]$, which is defined by $a\mapsto \varphi^*(a)$, 
$x\mapsto x$. Therefore $\psi:X_B\ra X$ is a purely inseparable 
finite morphism of degree $p^n$.

Since $f|_R:R\ra C$ is purely inseparable of degree $p^n$, so is 
$f|_R\circ\nu:R^\nu\ra C$, where $\nu:R^\nu\ra R$ is the normalization of $R$. 
We can also assume that $f|_R\circ\nu:R^\nu=C_{p^n}\ra C$ is the iterated 
$k$-linear Frobenius map. Assume $R=\Spec A_1$. Then $A$ is the integral 
closure of $A_1$ in its quotient field. Since $(f|_R\circ\nu)^*: A\ra A$ 
is the iterated $k$-linear Frobenius map, it is obviously injective. Hence 
$(f|_R)^*: A\ra A_1$ is also injective, which implies $A^{p^n}\subset A_1
\subset A$. Thus we have $R_B=\Spec(A\otimes_A A_1)$, where the $A$-module 
structures of $A$ and $A_1$ are induced by the iterated $k$-linear Frobenius 
map from $A$, and $g|_{R_B}:R_B\ra B$ is induced by $A\ra A\otimes_A A_1$, 
$a\mapsto a\otimes 1$. Define the morphism $s:B\ra R_B$ by $A\otimes_A A_1
\ra A$, $a\otimes b\mapsto ab$. We can verify that the morphism $s$ is 
well-defined and induces a section of $g:X_B\ra B$.
\end{proof}

The following lemma has already been proved in \cite[Proposition 3.3]{sc01}. 

\begin{lem}\label{4.8}
Let $C$ be a smooth projective curve, $f:X\ra C$ a $\PP^1$-bundle, 
and $R$ an irreducible curve on $X$. If $R^2<0$, then $f|_R:R\ra C$ 
is purely inseparable.
\end{lem}

\begin{proof}
Let $K$ be the separable closure of $K(C)$ in the field $K(R)$. 
Then $K(R)/K$ is purely inseparable and $K/K(C)$ is separable. 
Since $K/K(C)$ is a primitive extension, we can assume that 
$K=K(C)[x]/(h(x))$, where $h(x)$ is a monic irreducible separable 
polynomial over $K(C)$. Let $K'$ be the splitting field of $h(x)$. 
Then $K'/K(C)$ is a Galois extension. Note that there exists a 
surjective morphism $\varphi:B\ra C$ from a smooth projective 
curve $B$ such that $K(B)=K'$. Let $\xi$ be the generic point of 
the curve $R$, $\psi^{-1}(\xi)=\{\eta_1,\cdots,\eta_n\}$ the generic 
points of the preimage $\psi^{-1}(R)$, and $R_i$ the irreducible 
curves on $X_B$ determined by $\eta_i$. It is easy to see that 
$K(B)\otimes_{K(C)} K$ is a product of copies of $K(B)$, 
and the factors are just the function fields of the curves $R_i$. 
Therefore the set $\{\eta_1,\cdots,\eta_n\}$ is permuted by the 
Galois group $\Gal(K(B)/K(C))$ and $\deg h(x)=n$. Hence we have 
$R_1^2=\cdots=R_n^2$. If the self-intersection is nonnegative, 
then we have $(R_1+\cdots+R_n)^2\geq 0$, which contradicts $R^2<0$. 
Thus $R_1^2=\cdots=R_n^2<0$. Since $\rho(X_B)=2$, we have 
$\deg h(x)=n=1$, hence $K=K(C)$ and $K(R)/K(C)$ is purely inseparable.
\end{proof}

Let us start the proof of Theorem \ref{4.1}. First of all, we need 
some reductions.

Let $B=\sum b_iB_i$. 
We may assume that each irreducible component of $\Supp(B)$ is 
of negative self-intersection. Indeed, if there is a component $B_0$ 
of $B$ such that $B_0^2\geq 0$, then we define $B'=B-b_0B_0$ and 
$H'=H+b_0B_0$. Note that $D=K_X+B+H=K_X+B'+H'$, $(X,B')$ is KLT and 
$H'=D-(K_X+B')$ is ample since $B_0$ is nef. Therefore we have only to 
prove Theorem \ref{4.1} for the triple $(X,B',D)$.

If $B=0$ then $H=D-K_X$ is an ample integral divisor, hence we have 
$H^1(X,D)=H^1(X,-H)^{\vee}=0$ by Lemma \ref{4.2}, which contradicts 
the assumption $H^1(X,D)\neq 0$. Therefore we have $B\neq 0$.

Since $\rho(X)=2$, there exists at most one irreducible curve of 
negative self-intersection on $X$. This fact will be used frequently 
in the later argument. Therefore, $\Supp(B)$ consists of just one 
irreducible curve of negative self-intersection. We fix the following 
notation throughout the proof of Theorem \ref{4.1}:
\vskip 2mm
{\it $B=cC'$, where $0<c<1$ and $C'$ is an irreducible curve on $X$ 
with $C^{'2}<0$.}
\vskip 2mm
Let $E$ be a section of $f: X\ra C$, $\EE=f_*\OO_X(E)$ and 
$\LL=f_*\OO_E(E)$. Then $\EE$ is a rank two locally free sheaf on $C$ 
such that $X=\PP(\EE)$ and $\LL$ is an invertible sheaf on $C$, 
which fit into the following exact sequence:
\begin{eqnarray}
0\ra \OO_C\ra \EE\ra \LL\ra 0. \label{es:7}
\end{eqnarray}
Denote $-e=E^2=\deg\LL=\deg\EE$, $F$ the fibre of $f$. 
We divide into two cases by $e$.

\begin{case}[A]
Assume $e\geq 0$.
\end{case}

We shall prove that Case (A) cannot occur, namely, $f:X\ra C$ has no 
section of self-intersection less than or equal to 0.

Assume $C'\equiv aE+bF$. If $C'\neq E$, then we have $a\geq 0$ and 
$b\geq ae\geq 0$, hence $C^{'2}=a(2b-ae)\geq 0$, which is absurd. 
Therefore we have $C'=E$ and $-e=E^2=C^{'2}<0$. Assume $H\equiv xE+yF$. 
Since $H$ is ample, we have $x>0$ and $y>xe>0$.

By Serre duality, we have $h^1(F,D|_F)=h^0(F,K_F-D|_F)=h^0(F,(K_X-D)|_F)
=h^0(F,-(H+B)|_F)=0$. Hence $R^1f_*\OO_X(D)=0$, $f_*\OO_X(D)$ is a locally 
free sheaf on $C$ by Grauert's theorem. It follows from the Leray spectral 
sequence that $H^1(X,D)=H^1(C,f_*\OO_X(D))$. 
By Serre duality, we have $H^1(C,f_*\OO_X(D))=H^0(C,(f_*\OO_X(D))^\vee
\otimes\omega_C)^\vee=H^0(C,(f_*\OO_X(D-f^*K_C))^\vee)^\vee$. 

Note that $D-f^*K_C=K_{X/C}+H+B\equiv (x+c-2)E+(y-e)F$ and $f:X\ra C$ 
is a $\PP^1$-bundle, so there exists a divisor $G$ on $C$ of degree 
$y-e$ such that $D-f^*K_C=(x+c-2)E+f^*G$. 
If $x+c-2<0$ then we have $f_*\OO_X(D-f^*K_C)=0$, 
hence $H^1(X,D)=H^0(C,(f_*\OO_X(D-f^*K_C))^\vee)^\vee=0$, which 
contradicts the assumption $H^1(X,D)\neq 0$. Thus we have $x+c-2\geq 0$.

Note that $H^0(C,(f_*\OO_X(D-f^*K_C))^\vee)=
H^0(C,S^{(x+c-2)}(\EE)^\vee\otimes\OO_C(-G))\neq 0$, 
so a nonzero section $s$ determines the following exact sequence:
\[ 0\ra \OO_C\stackrel{s}{\ra}\cS\ra \cQ\ra 0, \]
where $\cS=S^{(x+c-2)}(\EE)^\vee\otimes\OO_C(-G)$ and $\cQ$ is the 
cokernel of $s$. Consider the inverse image in $\cS$ of the torsion 
subsheaf of $\cQ$, we can denote it by $\OO_C(G_1)$, where $G_1$ is an 
effective divisor on $C$ of degree $d_1\geq 0$. Denote the quotient 
$\cS/\OO_C(G_1)$ by $\cQ_1$, then $\cQ_1$ is a locally free sheaf on $C$ 
and we have the following exact sequence:
\[ 0\ra \OO_C(G_1)\ra \cS\ra \cQ_1\ra 0. \]
Taking the dual of this sequence, we have a surjective homomorphism 
$\cS^\vee\ra\OO_C(-G_1)$, which gives rise to the exact sequence 
$S^{(x+c-2)}(\EE)\ra \OO_C(-G-G_1)\ra 0$.

Applying Lemma \ref{4.3} to $S^{(x+c-2)}(\EE)\ra \OO_C(-G-G_1)\ra 0$, 
we have the inequality $-(y-e)-d_1\geq -(x+c-2)e$, i.e. 
$-d_1\geq (y-xe)+(1-c)e>0$, which contradicts $d_1\geq 0$. 
Therefore Case (A) cannot occur.

\begin{case}[B]
Assume $e<0$.
\end{case}

From the former argument, it follows that $f:X\ra C$ has no section of 
self-intersection less than or equal to 0. So we can assume that 
$X=\PP(\EE)$ with $\EE$ normalized and there is a section $E$ such that 
$\OO_X(E)=\OO_X(1)$ and $E^2=\deg\EE=-e>0$. 
Therefore we have the following exact sequence:
\begin{eqnarray}
0\ra \OO_C\ra \EE\ra \LL\ra 0, \label{es:8}
\end{eqnarray}
where $\LL$ is an ample line bundle on $C$.

By Lemma \ref{4.8}, we have $f|_{C'}: C'\ra C$ is purely inseparable 
of degree $p^n$. If $n=0$ then $C'$ is a section of $f: X\ra C$ with 
negative self-intersection, which is absurd by the former argument. 
Therefore $f|_{C'}: C'\ra C$ is the iterated $k$-linear Frobenius map 
of degree $p^n>1$. Furthermore, the exact sequence (\ref{es:8}) 
is not split. Otherwise, there would be a section of $f: X\ra C$ with 
negative self-intersection, which is absurd by the former argument.

\begin{subcase}[B--1]
Assume that $E$ is not ample.
\end{subcase}

In this subcase, since $E$ is nef and big, but not ample, we have an 
irreducible curve $C_1$ such that $E.C_1=0$, hence $C_1^2<0$ by the 
Hodge index theorem. Thus we have $C_1=C'$ and $E\cap C'=\emptyset$. 
By Lemma \ref{4.4} (ii) and (iii), $E_B$ and $C'_B$ are disjoint 
sections of the $\PP^1$-bundle $g:X_B\ra B$, where $\varphi:B\ra C$ 
is the iterated $k$-linear Frobenius map of degree $p^n$. 
This implies that after pulling back by $\varphi:B\ra C$, 
the exact sequence (\ref{es:8}) becomes split 
(cf.\ \cite[Exercise V.2.2]{ha}). Therefore, we can assume that 
there exists an $m\geq 1$ such that the pulling back of (\ref{es:8}) 
by $F^{m-1}$ is not split, while the pulling back of (\ref{es:8}) 
by $F^{m}$ is split.

Tensoring the exact sequence (\ref{es:1}) by $\LL^{-p^{m-1}}$ and taking 
cohomology groups, we have the following exact sequence:
\[ 0\ra H^0(C,\BB^1\otimes\LL^{-p^{m-1}})\ra H^1(C,\LL^{-p^{m-1}})\stackrel{F^*}\ra H^1(C,\LL^{-p^m}). \]
Let $\alpha\in H^1(C,\LL^{-1})$ be the extension class of the 
exact sequence (\ref{es:8}). By assumption, we have $0\neq F^{m-1*}(\alpha)
\in H^1(C,\LL^{-p^{m-1}})$ and $0=F^{m*}(\alpha)\in H^1(C,\LL^{-p^m})$. 
Therefore, we have 
$0\neq F^{m-1*}(\alpha)\in H^0(C,\BB^1\otimes\LL^{-p^{m-1}})$ 
by the above exact sequence. 
As a consequence, we have $H^0(C,\BB^1\otimes\LL^{-p^{m-1}})\neq 0$, 
hence $n(C)>0$ by Lemma \ref{2.5}.

\begin{subcase}[B--2]
Assume that $E$ is ample.
\end{subcase}

Since $\EE$ is normalized and $E$ is ample, it is easy to see that 
all of sections of $f:X\ra C$ are ample, which completes the proof 
of Theorem \ref{4.1}.
\vskip 3mm

We give some remarks on Subcase (B--2).

\begin{rem}\label{4.5}
(i) Subcase (B--2) occurs only if $g(C)\geq 2$. Indeed, if $C\cong \PP^1$ 
then we would have $e\geq 0$ by \cite[Corollary V.2.13]{ha}, which contradicts 
the assumption $e<0$. If $g(C)=1$ and write $C'\equiv p^nE+bF$, then we have 
$b\geq \frac{1}{2}p^ne$ by \cite[Proposition V.2.21]{ha}. Hence we would 
have $C^{\prime 2}\geq 0$, which contradicts the assumption $C^{\prime 2}<0$.

(ii) The pulling back of the exact sequence (\ref{es:8}) by $F^m$ can never 
be split for any $m\in\N$. Indeed, if we assume the contrary then $\OO_C$ 
would be a direct summand of the ample vector bundle $F^{m*}\EE$, hence 
would be ample, which is absurd.

(iii) By \cite[Theorem 3.5]{sc01}, we can show that 
there is an irreducible curve $E'$ on $X$ such that $f|_{E'}:E'\ra C$ 
is purely inseparable of degree $p^m$, $E^{\prime 2}>0$ and $E'\cap C'
=\emptyset$. Thus $E'$ is not a section of $f$. If $m<n$ then we can 
follow a similar argument to Subcase (B--1) to show $n(C)>0$. 
However, the case $m\geq n$ is complicated, and it seems that 
we cannot obtain any contradiction from the existence of such $E'$.
\end{rem}

The following corollary is a consequence of the proof of Theorem \ref{4.1}.

\begin{cor}\label{4.6}
Let $f:X\ra C$ be a $\PP^1$-bundle over a smooth projective curve $C$. 
Then the Kawamata-Viehweg vanishing theorem holds on $X$ if one of the 
following conditions holds:

(i) One of sections of $f$ has self-intersection less than or equal to 0,

(ii) $n(C)\leq 0$ and one of sections of $f$ is not ample,

(iii) $g(C)\leq 1$.
\end{cor}

Note that in Subcase (B--2), the ampleness of $E$ is equivalent 
to the ampleness of $\EE$ by \cite[Proposition 3.2]{ha66}, 
so we may take the following problem into account.

\begin{prob}\label{4.7}
Let $f:X=\PP(\EE)\ra C$ be a $\PP^1$-bundle, where $\EE$ is a rank two 
locally free sheaf on $C$ with $g(C)\geq 2$. Assume that $\EE$ is normalized 
and ample. Does the Kawamata-Viehweg vanishing theorem hold on $X$?
\end{prob}

Although there is no evidence, we expect that Problem \ref{4.7} 
holds true. Problem \ref{4.7}, together with Theorem \ref{4.1}, 
implies the following conjecture, which would give a 
characterization of counterexamples to the Kawamata-Viehweg 
vanishing theorem on geometrically ruled surfaces.

\begin{conj}\label{4.9}
Let $C$ be a smooth projective curve. Then there exists a counterexample 
to the Kawamata-Viehweg vanishing theorem on some geometrically ruled 
surface over $C$ if and only if $C$ is a Tango curve.
\end{conj}

Note that Conjecture \ref{4.9} is a special case of Conjecture \ref{1.2}. 
On the other hand, it is of certain interest to consider the following 
conjecture, which is the converse of Theorem \ref{3.5} and similar 
to Conjecture \ref{1.2}.

\begin{conj}\label{4.10}
If there is a counterexample to the Kodaira vanishing theorem on a smooth 
projective surface $Y$, then there exists a fibration $g:Y\ra C$ to a 
smooth projective curve $C$ such that the general fibre of $g$ is a singular 
rational curve and $C$ is a Tango curve of integral type.
\end{conj}

\textsc{School of Mathematical Sciences, Fudan University, 
Shanghai 200433, China}

\textit{E-mail address}: \texttt{qhxie@fudan.edu.cn}

\end{document}